\documentclass[a4paper,12pt,reqno]{amsart}
\usepackage{amssymb,amsfonts,amscd,latexsym}
\usepackage{xy}
\xyoption{matrix}
\xyoption{arrow}

\newcommand{\ol}{\overline}
\newcommand{\rest}[2]{{{#1}_{\kern-.5pt|{#2}}}}

\newcommand{\N}{\mathbb{N}}

\newcommand{\Q}{\mathbb{Q}}
\newcommand{\mloc}{M_{\text{\rm loc}}(A)}
\newcommand{\Mloc}[1]{M_{\text{\rm loc}}({#1})}
\newcommand{\eps}{\varepsilon}
\newcommand{\Prim}{\text{Prim}(A)}
\newcommand{\longrightarrowraised}{{\hbox{\raise.5\jot\hbox{\scriptsize$\longrightarrow$}}}}
\newcommand{\dirlim}{{\smash{\underset{\longrightarrowraised}
                   {\operatorname{lim}}}\vphantom{a^{}_{a_f^{}}}}\,}

\def\C*{{\sl C*}-algebra}
\def\Cs*{{\sl C*}-subalgebra}
\def\AF/{{\sl A\kern-.5pt F}-algebra}
\def\AW*{{\sl AW*}-algebra}

\newtheorem{lem}{Lemma}[section]

\newtheorem{theor}[lem]{Theorem}
\newtheorem{prop}[lem]{Proposition}

\theoremstyle{definition}

\newtheorem{rema}[lem]{Remark}

\begin{document}

\title{A not so simple local multiplier algebra}
\author{Pere Ara and Martin Mathieu}
\address{Departament de Matem\`atiques, Universitat Aut\`onoma de
Barcelona, 08193 Bellaterra (Barcelona), Spain}
\email{para@mat.uab.es}
\address{Department of Pure Mathematics, Queen's University Belfast,
Belfast BT7 1NN, Northern Ireland}
\email{m.m@qub.ac.uk}

\thanks{The first-named author was partially supported by
DGI and European Regional Development Fund, jointly, through
Project BFM2002-01390, and by the Comissionat per Universitats i
Recerca de la Generalitat de Catalunya.
This paper is part of a research project supported by the Royal Society.}

\dedicatory{To the memory of Gert Kj\ae rgaard Pedersen.}


\begin{abstract}
We construct an \AF/ $A$ such that its local multiplier algebra
$\mloc$ does not agree with $M_{\text{loc}}(\mloc)$, thus answering
a question raised by G.~K.~Pedersen in~1978.
\end{abstract}

\maketitle

\baselineskip=15pt

\section{Introduction}\label{sec:intro}

For a \C* $A$ denote by $M(A)$ its multiplier algebra. A closed (two-sided) ideal
$I$ of $A$ is called {\it essential\/} if it has non-zero intersection with each
non-zero ideal of~$A$. Let $I,\,J$ be closed essential ideals of $A$ such that $J\subseteq I$.
Then the restriction mapping induces a *-monomorphism $M(I)\to M(J)$. The direct limit
of all $M(I)$ along the downward directed family of all closed essential ideals
and with these connecting mappings is the local multiplier algebra $\mloc$
of $A$, first introduced by Pedersen in~\cite{ped}.
Further properties of $\mloc$ were studied in~\cite{AM2}.

Among the questions which were left open in \cite{ped} are the following.
\begin{enumerate}
\item Is every derivation of $\mloc$ inner?
\item Does the equality $\Mloc\mloc=\mloc$ hold?
\end{enumerate}
Pedersen showed that a derivation $d$ on $A$ can be (uniquely) extended to a derivation on
$\mloc$, and becomes inner in $\mloc$ provided $A$ is separable. (For a detailed
account of his argument, and related questions, see \cite[Section~4.2]{AM2}.)
Despite some interesting contributions by Somerset~\cite{somer}, Question~(1)
remains open. Note that a positive answer includes the classical results for
simple \C*s, for von Neumann algebras and for \AW*s by Sakai, Kadison and Olesen,
respectively (compare \cite[Section~8.6]{pedbook}).

If the answer to Question~(2) were positive, to prove~(1) it would suffice
to show that every derivation on $A$ becomes inner in $\mloc$.
This occurs in particular when $M(A)$ is an \AW* or $A$ is simple; for, in this case,
$\Mloc{M(A)}=M(A)$ and $\Mloc{M(A)}$ and $\mloc$ always coincide 
\cite[Section~2.3]{AM2}.
It also occurs when $\mloc$ is an \AW* or is simple; the former holds for
every commutative \C* \cite[3.1.5]{AM2} and the latter is indeed 
possible in non-trivial cases as was shown in~\cite{AM1}.
Further evidence for a positive answer is provided by the local Dauns--Hofmann
theorem which implies that $Z(\Mloc\mloc)=Z(\mloc)$ for every \C*~$A$
\cite[3.2.6]{AM2}. Somerset showed in \cite[Theorem~2.7]{somer}
that (2) holds for every unital separable \C* $A$ such that its primitive
spectrum $\text{Prim}(A)$ contains a dense $G_\delta$ of closed points;
hence in particular if $\text{Prim}(A)$ is Hausdorff.
Argerami and Farenick recently derived (2) under the assumption that $A$ is separable
and contains a minimal essential ideal of compact elements; in this case
$\mloc$ coincides with the injective envelope of $A$ and is a type~I von
Neumann algebra~\cite{AF}.

In general, however, it turns out that the answer to Question~(2) is negative.
In this paper we provide a class of examples to this effect. Our main result is
the following.
\begin{theor}\label{thm:main}
There exist unital, primitive \AF/s\/ $A$ such that\/ $\Mloc{\mloc}\ne\mloc$.
\end{theor}
The strategy to obtain such \AF/s follows the ideas in~\cite{AM1},
where we gave examples of non-simple \AF/s $A$ with the property that $\mloc$ is simple.
To specify an \AF/ it is, of course, enough to write down its {\sl K}-theoretic
invariant. It emerges, however, that working with the monoid $V(A)$ of equivalence
classes of projections in $M_\infty(A)$ gives us a better control on the 
order-theoretic properties. Since $V(A)$ is cancellative in this case,
this approach is of course equivalent to the usual {\sl K}-theoretic one;
however, it allows for a description of the ideal structure of the multiplier
algebras of closed essential ideals of~$A$ (which is the decisive step in understanding
$\mloc$). By work of Goodearl~\cite{Good} and Perera~\cite{perera},
for a $\sigma$-unital \C* $A$ of real rank zero and stable rank one, 
the monoid $V(M(A))$ can be completely described by the monoid of countably generated
complemented intervals on $V(A)$. In order to obtain a like description of
$V(\mloc)$, a localisation procedure is needed, which was carried out in
\cite[Theorem~2]{AM1}.
\begin{theor}\label{thm:from-old}
Let\/ $A$ be a unital \AF/. Then\/ $\mloc$ has real rank~zero
and\/ $V\bigl(\mloc\bigr)$ is isomorphic to\/ 
$\Lambda_{\text{\rm loc}}\bigl(V(A),[1_A]\bigr)$, the monoid of local intervals.
\end{theor}
All the necessary concepts and notation will be introduced in 
Section~\ref{sec:monoids}, where we shall construct a certain countable, abelian
monoid $M$ which, endowed with the algebraic order, leads to a localised monoid $M'$
(representing $V(\mloc)$) that has a unique minimal order-ideal.
As a result, $\mloc$ has a unique minimal closed ideal~$I$ so that
$\Mloc\mloc=M(I)$.

However, the tools available in the literature are not sufficient to
compute the structure of the projections in the \C*~$M(I)$. The
reason is that $I$ is not $\sigma$-unital and, moreover, the
projections in $I$ can fail to satisfy cancellation. To resolve
this problem we use a different technique in Section~\ref{sec:multi-algebras}, 
which is inspired by the
geometry of our examples. We construct a {\it sequence\/} of
projections in $\mloc$ strictly converging in the
$I$-topology to a projection in $M(I)\setminus \mloc$. This allows us to conclude
that $\Mloc{\mloc}\ne\mloc $.

Both parts of the construction of our example are fairly technical.
Thus, in a brief Section~4, we reflect on the nature of the
example and possible further studies on the ideal structure of the local
multiplier algebra.

\section{Monoids}\label{sec:monoids}

This section is devoted to the construction of an ordered monoid with very special
properties. These will be exploited when it comes to exhibit the structure
of the local multiplier algebra associated with the corresponding \AF/ in
the following section.

We begin by fixing our setting. Let $X$ be an infinite compact
metrizable space, and let $t_0$ be a non-isolated point in~$X$.
Denote by $C(X)$ the set of all continuous {\it real-valued\/}
functions on~$X$ equipped with pointwise order. Let $G$ be a countable,
additive subgroup of $C(X)$ with the following properties:

\begin{itemize}
\item[(i)] $G$ is a sublattice of $C(X)$ and $\,\Q\,G\subseteq G$;
\smallskip
\item[(ii)] $G$ contains a function $f_0$ such that $0\le f_0\le 1$,
            $f_0(t_0)=0$ and $f_0(t)>0$ for all $t\in X\setminus\{t_0\}$;
\smallskip
\item[(iii)] for each $f\in G$ there are an open neighbourhood $V$ of $t_0$
             and $\lambda,\,\mu\in\Q$ such that $f=\lambda+\mu f_0$ on~$V$;
\smallskip
\item[(iv)] for every closed subset $K\subseteq X$, every open subset $V$ with
            $K\subseteq V$ and every $\rho\in\Q_+$, there are an open subset $U$
            with $K\subseteq U\subseteq V$ and $r\in G$ such that
            $0\le r\le\rho$ on~$X$, $r=\rho$ on $U$ and $r=0$ on ${^c}V=X\setminus V$.
\end{itemize}

\noindent
Property~(iv) requires $G$ to contain enough ``Urysohn functions".
It implies in particular that $1\in G$ (take $K=X$ and $\rho=1$).

\begin{prop}\label{group}
There exists a countable subgroup\/ $G$ of\/ $C(X)$ with the above properties\/
{\rm(i)--(iv)}.
\end{prop}

\begin{proof}
We can take a countable set $T$ of Urysohn functions such that,
for every $f\in T$, either $f=0$ or $f=1$ on a neighbourhood of
$t_0$ as follows. For each $n\in\N$, take open balls
$U^{(n)}_0,U^{(n)}_1, \dots ,U^{(n)}_{k_n}$ of radius $1/n$ and
centres $t^{(n)}_i$, for $i=1,\dots ,k_n$, with $t^{(n)}_0=t_0$,
which form a cover of~$X$. Then consider pairs of open subsets
$U$ and $V$ such that $\ol{V}\subseteq U$, and such that $V$ and
$U$ are finite unions of some of the open balls $U^{(m_j)}_j$,
but only those pairs $(U,V)$ such that either
$t_0\in V$ or $t_0$ belongs to the interior of $X\setminus U$. For
each such pair $(U,V)$, choose a Urysohn function
$f_{(U,V)}\colon X\to [0,1]$ such that $f_{(U,V)}=1$ on $V$ and
$f_{(U,V)}=0$ on $X\setminus U$. The set $T$ is defined as the set
of all these functions~$f_{(U,V)}$. It is clear that each function
in $T$ is either $0$ or $1$ on a neighbourhood of $t_0$, and it is
not hard to see that the set $\{\mu f\mid f\in T,\,\mu\in\Q_+\}$
contains enough Urysohn functions in the sense of property~(iv).

Set $G_1=\Q\,T+\Q\,f_0+\Q\,1$, the $\Q$-linear span of $T$, $f_0$ and
$1$, and observe that, for $f\in G_1$, there are rational numbers
$\lambda$ and $\mu$ such that $f=\lambda +\mu f_0$ on a
neighbourhood of~$t_0$.

Let $T_1$ be the set consisting of all the functions of the form
$f\wedge g$ and $f\vee g$, for $f,\,g\in G_1$. It is then clear that
each function in $T_1$ is of the form $\lambda +\mu f_0$, for some
$\lambda,\,\mu\in\Q$, on a neighbourhood of~$t_0$.

Proceeding inductively, suppose that we have a countable set of
functions $T_n$ with the property that, for each $f\in T_n$, there
are $\lambda,\mu \in \Q$ such that $f=\lambda +\mu f_0$ on a
neighbourhood of~$t_0$. Then define $G_{n+1}$ as the $\Q$-linear
span of $T_n$, and define $T_{n+1}$ as the set of all functions
$f\wedge g$, $f\vee g$, for $f,g\in G_{n+1} $. Clearly
$T_{n+1}$ enjoys the same property as~$T_n$.  Finally, set
$G=\bigcup_{n=1}^{\infty}G_n$. Then $G$ satisfies the desired conditions.
\end{proof}

From now on, $G$ will denote a countable subgroup of $C(X)$ satisfying
the above properties~(i)--(iv).
Whenever $t\in X$, we write $V(t)$ to denote {\it some\/} open neighbourhood
of~$t$ and $V[t]$ to denote the punctured neighbourhood $V(t)\setminus\{t\}$.
Let
\[
M=\{f\in G\mid f\geq0,\;f(t)>0\,\text{ on }V[t_0]\}\cup\{0\}.
\]
Note that $M$ is a countable, additive monoid, closed under
multiplication by positive rational numbers. The {\it algebraic
order\/} in $M$ will be denoted by $\le_M$. We fix a canonical
order-unit $u=1$ in~$M$.

We recall that the {\it algebraic pre-order\/} on an abelian monoid $L$
is defined by
\[
x\leq_L y\,\text{ if }\,y=x+z\,\text{ for some }\,z\in L.
\]
In the case where $L$ is cancellative, $\leq_L$ is a (partial) order. We write
$x<_L y$ if $x\leq_L y$ and $x\ne y$. We also recall a few order-theoretic concepts
that will be used in the following. (For more details, see~\cite{AM1} and~\cite{perera}.)

An {\it order-ideal\/} of $(L,\leq_L)$ is a hereditary submonoid;
an {\it order-unit\/} is an element such that $L$ is the smallest order-ideal
containing it; an {\it interval\/} is an upward directed hereditary non-empty
subset of~$L$. The monoid $L$ is said to be {\it prime\/} if each pair of
non-zero order-ideals of $L$ has non-zero intersection. Suppose that $L$ is
a {\it Riesz monoid}; that is, whenever $x,y_1,y_2\in L$ satisfy
$x\leq_L y_1+y_2$ there exist $x_1,x_2\in L$ such that $x=x_1+x_2$ and $x_i\leq_L y_i$
for $i=1,2$. Then the sum $E+F$ of two intervals $E$ and $F$ in $L$ is defined by
\[
E+F=\{x+y\mid x\in E,\,y\in F\}
\]
and is an interval in~$L$. Let $D$ be a fixed interval. The interval $E$ is said to be
{\it complemented\/} ({\it with respect to}~$D$) if there are an interval $F$ and
some $k\geq1$ such that $E+F=kD$. We denote by
$\Lambda(L,D)$ the abelian monoid of all complemented intervals in~$L$.

The following important localisation procedure will be applied to various
Riesz mono\-ids in the sequel. Suppose $L$ is a prime Riesz monoid with order-unit~$v$.
Let $N$ be an order-ideal in~$L$; then its {\it canonical interval\/} $D_N$
is defined by $D_N=\{x\in N\mid x\leq_L v\}$. Let $N_1, N_2$ be order-ideals
in $L$ with $N_1\subseteq N_2$. The restriction mapping
\[
\phi_{N_1,N_2}\colon\Lambda(N_2,D_{N_2})\longrightarrow\Lambda(N_1,D_{N_1})
\]
defined by
$\phi_{N_1,N_2}(E)=E\cap N_1$ is a monoid homomorphism.
Whenever \hbox{$N_1\subseteq N_2\subseteq N_3$} is a chain of order-ideals, we have
$\phi_{N_1,N_3}=\phi_{N_1,N_2}\phi_{N_2,N_3}$.
Therefore we can define the {\it monoid of local intervals\/}
$\Lambda_{\text{\rm loc}}(L,v)$ of $(L,v)$ as the direct limit
of the family
\[
\bigl\{\Lambda(N,D_N);\,\phi_{N_1,N_2},\,N_1\subseteq N_2\bigr\},
\]
where $N$ runs through the downward directed set of all non-zero
order-ideals of~$L$.

This procedure will now be applied to the monoid $(M,u)$.

The proof of our first lemma is exactly the same as the one of the
first part of Theorem~3 in~\cite{AM1} and hence is omitted.

\begin{lem}\label{Riesz} 
The monoid\/ $(M,\leq_M)$ is a prime cancellative Riesz monoid.
\end{lem}

For each non-zero $f\in M$, in order to simplify the notation, set
\begin{equation*}
\begin{split}
N_f   &=\{g\in M\mid g\le _M nf \text{ for some } n\ge 1\},\\
N_f^* &=N_f\setminus\{0\},\\
D_f   &=\{g\in N_f\mid g\le _M u\},\\
L_f   &=\{g\in C_b(U_f)_+\mid\exists\,z\in N_f^*\colon z\leq g\,\text{ on }U_f\}\cup\{0\},
\end{split}
\end{equation*}
where $U_f$ stands for the co-zero set of~$f$. For $f'\le _M f$ we have a canonical map
\[
\phi_{f',f}\colon \Lambda(N_f,D_f)\to \Lambda (N_{f'},D_{f'})
\]
defined by $\phi_{f',f}(E)= E\cap N_{f'}$. (That is,
$\phi_{f',f}=\phi_{N_{f'},N_f}$.)

\smallskip
We will denote by $M'$ the prime Riesz monoid
$M'=\Lambda_{\text{loc}}(M,u)=\dirlim\Lambda(N_f,D_f)$; in general, this may not be
cancellative. For an interval $E$ in $\Lambda(N_f,D_f)$, where $f\in M\setminus\{0\}$,
we denote the class of $E$ in $M'$ by~$\overline{E}$. There is a canonical order-unit
in $M'$ given by $u'=\overline{[0,u]}$.
Let $J$ be the order-ideal of $M'$ generated by $\overline{[0,f_0]}$.

One of the key properties of $G$, as stated in Proposition~\ref{group},
is that each $f$ in $G$ is of the form $\lambda +\mu f_0$ on $V(t_0)$
for some $\lambda,\,\mu\in\Q$. This obviously implies that, given $f\in M$
such that $f(t_0)=0$,
there is a rational number~$\mu$ such that $f=\mu f_0$ on~$V(t_0)$.

\begin{prop}\label{prop:min-order-ideal}
The monoid\/ $M'$ has a unique minimal order-ideal\/ $J$, the order-ideal
generated by\/ $\ol{[0,f_0]}$.
\end{prop}

\begin{proof}
It suffices to show that, for every non-zero $x\in M'$, we have
$\ol{[0,f_0]}\le_{M'}nx$ for some $n\in\N$. Let
$E\in\Lambda(N_f,D_f)$ be a representative of $x$, where $0\ne f\le _M f_0$.
Take a non-zero element $g$ in $E$. Since
$g(t_0)=0$, there is a rational number $\mu>0$ such that $g=\mu f_0$
on $V(t_0)$. Thus, on $V[t_0]$, we have $f_0 \ll ng$ for some
$n\in\N$. Take $f'\in M\setminus \{0\}$ with $f'\le _M f$ such
that $U_{f'}\subseteq V(t_0)$. Observe that
\[
[0,f_0]\cap N_{f'}+[0,(n\mu-1) f_0]\cap N_{f'}=[0,n\mu f_0]\cap
N_{f'}=[0,ng]\cap N_{f'}.
\]
On the other hand, $ng\in nE$ so that $nE=[0,ng]+E'$, where $E'$ is
the interval in $M$ defined as
\[
E'=\{z\in M\mid z+ng\in E\}.
\]
It follows that $[0,ng]\cap N_{f'}+E'\cap N_{f'}=nE\cap N_{f'}$,
and so
\[
[0,f_0]\cap N_{f'}+
\bigl([0,(n\mu-1)f_0]\cap N_{f'}+E'\cap N_{f'}\bigr)=n(E\cap N_{f'}),
\]
which shows that $\ol{[0,f_0]}\le_{M'}n\ol{E}=nx$, as desired.
\end{proof}

In Section~\ref{sec:multi-algebras} we shall need a functional representation
of the monoid~$M'$. Let  $f$ be
a non-zero element of $M$ such that $f\le _M f_0$. Note that the
set of order-ideals $N_{f}$ with such $f$ is cofinal;
so in order to study $M'$, we may restrict attention to those
$\Lambda(N_f,D_f)$. Fix such an element~$f$. Then, for every $z\in N_f$,
we have $z(t_0)=0$, so our hypothesis gives that for some
$\mu\in\Q^+$ we have $z=\mu f_0$ on~$V(t_0)$.

For $h\in C_b(U_f)_+$ we set
\[
I_f(h)=\{g\in N_f\mid\exists\,z\in N_f\colon g<_M z,\ z\ll h\text{ on }
V[t_0],\ z\leq h\text{ on }U_f\} \cup \{0\}.
\]
The following description of $I_f(h)$ will be used subsequently
several times without specific reference.

\begin{lem}\label{lem:ifh-interval}
For every\/ $h\in L_f\setminus\{0\}$,
\begin{equation*}
\begin{split}
I_f(h) &=\{g\in N_f\mid g\leq h\text{ and }g\ll g'\ll h\text{ on }V[t_0]
\text{ for some }g'\in N_f\}\\
       &=\{g\in N_f\mid g<_M z\leq h\text{ for some }z\in N_f\}.
\end{split}
\end{equation*}
Moreover\/ $I_f(h)$ is an interval in\/ $N_f$.
\end{lem}

\begin{proof}
Put $\tilde I_f(h)=\{g\in N_f\mid g\leq h\text{ and }g\ll g'\ll h
\text{ on }V[t_0]\text{ for some }g'\in N_f\}$. It is evident that
$I_f(h)\subseteq \tilde I_f(h)$.

To show the reverse inclusion, assume that $g\in N_f$ is such that
$g\leq h$ and $g\ll g'\ll h$ on $V[t_0]$ for some $g'\in N_f$.
Take an open neighbourhood $V$ of $t_0$ with $\ol V\subseteq V(t_0)$.
Let $r\in M$ be such that $r=\rho\gg g'$ on $V$, $r=0$ on ${}^c V(t_0)$
and $0\leq r\leq\rho$ for some $\rho\in\Q$. Then
$g'\wedge r=g'$ on~$V$, $g'\wedge r\ll h$ on $V[t_0]$ and $g'\wedge r=0$
on ${}^cV(t_0)$. Let $z=(g'\wedge r)\vee g\in N_f$.
Then $g<_M z$, $z\ll h$ on $V[t_0]$ and $z\leq h$ on~$U_f$;
thus $g\in I_f(h)$.

Put $I'_f(h)=\{g\in N_f\mid g<_M z\leq h\,\text{ for some }z\in N_f\}$;
clearly $I_f(h)\subseteq I'_f(h)$. On the other hand, if
$g\in I'_f(h)$ and $z\in N_f$ satisfies $z\leq h$ on $U_f$ and
$g<_M z$ then we take $r\in M$ such that $r=\rho\gg z-g$ on~$V$,
an open neighbourhood of~$t_0$ with $\ol V\subseteq V(t_0)$ and
$g\ll z$ on $V[t_0]$, $r=0$ on ${}^cV(t_0)$ and $0\leq r\leq\rho$.
Upon replacing $z$ by $\frac12(z-g)\wedge r +g\in N_f$ we find
that $g\in I_f(h)$.

Clearly $I_f(h)$ is a non-empty hereditary subset of~$N_f$. To show that
it is upward directed take $g_1,\,g_2\in I_f(h)$.
There are $\mu_1,\,\mu_2\in\Q$ such that $g_i=\mu_if_0$ on $V(t_0)$,
$i=1,2$. We may assume that $\mu_1\geq\mu_2$. There exists $g'\in N_f$
such that $g_1<_Mg'$, $g'\ll h$ on $V'[t_0]$, where $V'(t_0)\subseteq V(t_0)$,
and $g'\leq h$ on~$U_f$.
Take an Urysohn function $r\in M$ as above, where $r=\rho\gg g'-g_1$
on~$V$ with $\ol V\subseteq V'(t_0)$ and $r=0$ on ${}^cV'(t_0)$. Set 
$g''=\frac12(g'-g_1)\wedge r+ g_1\vee g_2$.
Then $g_1\leq_M g''$, $g_2\leq_M g''$ and for
$z=\bigl((g'-g_1)\wedge r\bigr)+g_1\vee g_2\in N_f$ we have
$g''<_M z$, $z\ll h$ on $V'[t_0]$ and $z\leq h$ on~$U_f$. Thus $g''\in I_f(h)$.
\end{proof}

Under the standing assumption that $f\in M\setminus\{0\} $ with $f\le _M f_0$
is given, we will now define a monoid homomorphism
$\tau_f\colon\Lambda(N_f,D_f)\to L_f$. For $E\in\Lambda(N_f,D_f)$ let
$\tau_f(E)=\sup E$ be the pointwise supremum over all functions in~$E$,
restricted to~$U_f$. Then $\tau_f$ has the following properties.

\smallskip
(1)\quad$\tau_f(D_f)=1$ on $U_f$;

\smallskip\noindent
This follows easily from the existence of sufficiently many Urysohn
functions in~$G$.

\smallskip
(2)\quad$\tau_f(E_1+E_2)=\tau_f(E_1)+\tau_f(E_2)$ for all
$E_1,\,E_2\in\Lambda(N_f,D_f)$;

\smallskip\noindent
This is straightforward. As a consequence of~(1) and~(2), $\tau_f(E)$ is
a continuous function on~$U_f$ and $\tau_f$ is a monoid homomorphism.

\smallskip
(3)\quad$h=\sup I_f(h)$ for each $h\in L_f$.

\smallskip\noindent
It follows from the definition of $I_f(h)$ that $h\geq\sup I_f(h)$.
In order to show the converse inequality, suppose that
$h\ne0$, let $\eps>0$ and take $z\in N_f^*$ with $z\leq h$. Let
$t\in U_f$ be such that $h(t)>0$. Let $V$, $V_0$
be disjoint open neighbourhoods of $t$ and $t_0$, respectively,
with the property that $\ol V\subseteq U_f$,
$h(s)-\eps<\rho<h(s)$ for all $s\in V$ and some $\rho\in\Q_+$
and $z\gg0$ on~$V_0\setminus\{t_0\}$.
There is $r\in M$ such that $0\leq r\leq\rho$, $r=\rho$ on~$W$ and
$r=0$ on~${}^cV$, where $W$ is some open neighbourhood of $t$ with
$\ol W\subseteq V$.
There is $r'\in M$ such that $0\leq r'\leq1$, $r'=1$ on~$W'$ and
$r=0$ on~${}^cV_0$, where $W'$ is some open neighbourhood of $t_0$ with
$\ol{W'}\subseteq V_0$. Then
$g=\frac14(z\wedge r')\vee r\in N_f$. In fact, $g\in I_f(h)$ since,
for $g'=\frac12(z\wedge r')\in N_f$, we have $g\ll g'\ll h$ on $W'\setminus\{t_0\}$.
As $g(t)=r(t)=\rho>h(t)-\eps$ it follows that $h=\sup I_f(h)$.

\medskip
Let $f'\in M\setminus\{0\}$, $f'\leq_M f$. Let $R\colon L_f\to L_{f'}$
denote the restriction mapping. (Note that $h_{|U_{f'}}\in L_{f'}$ if $h\in L_f$.
For $h=0$ this is evident, so assume that $h$ is non-zero.
By definition, there is $z\in N_f^*$ with $z\leq h$ on~$U_f$
wherefore $z\wedge f'\in N_{f'}^*$ and $z\wedge f'\leq h$ on~$U_{f'}$.)

We obtain the following commutative diagram.
\begin{equation}
\begin{split}
\xymatrix{\ar[d]_{\phi_{f',f}} \Lambda(N_f,D_f)\ \;\ar[r]^{\quad\ \tau_f}
            &\ \; L_f\ \ar[d]^R\\
          \Lambda(N_{f'},D_{f'})\ \ar[r]_{\quad\ \tau_{f'}\vphantom{T}}
            &\ \; L_{f'}}
\end{split}
\tag{2.1}
\end{equation}

To verify the commutativity of the diagram take $E\in\Lambda(N_f,D_f)$
and note at first that
\[
\tau_{f'}\circ\phi_{f',f}(E)=\sup(E\cap N_{f'})\leq\sup E_{|U_{f'}}.
\]
For the reverse inequality, let $h=\sup E$ and let $t\in U_{f'}$ such that $h(t)>0$.
Let $\eps>0$. By Property~(3) above, there is $g'\in I_{f'}(h)$ such that
$h(t)\geq g'(t)>h(t)-\eps$. By definition of~$h$, there is $z\in E$ such that
$z(t)\geq g'(t)$. Put $g=g'\wedge z\in E\cap N_{f'}=\phi_{f',f}(E)$.
As $g(t)=g'(t)$ it follows that $\sup\phi_{f',f}(E)\geq h_{|U_{f'}}$.

Set $L=\dirlim L_f$, where the morphisms $L_f\to L_{f'}$ for $f'\le _M f$
are given by restriction. The commutativity of the above diagram enables
us to define a monoid homomorphism
\[
\tau\colon M'\to L
\]
via the maps $\tau_f\colon\Lambda_f(N_f,D_f)\to L_f$. In contrast to
\cite[Theorem~3]{AM1}, $\tau$ may not be injective in general and hence
$M'$ need not be cancellative.

In order to describe the image of $\tau$ we use a similar approach as in~\cite{AM1}.
For each non-zero $f\in M$, the function $h\in L_f$ has the property~(C) provided that,
for each $t\in U_f\setminus\{t_0\}$, there exists $z_t\in N_f$ such that $z_t\le h$ on $U_f$,
$h=z_t$ on $V(t)$ and $z_t\ll h$ on $V[t_0]$. The subset $C$ of $L$ consists of those
elements $y\in L$ such that there is a representative $h\in L_f$ of $y$ with property~(C),
together with $y=0$. Evidently $C$ is a submonoid of~$L$.

\goodbreak
\begin{lem}\label{oembedding}
Let\/ $x,y\in C$. If\/ $x<_L y\,$ then\/ $y-x\in C$.
\end{lem}

\begin{proof}
Without restricting the generality we can assume that both $x$ and $y$ are non-zero.
By hypothesis, there is $x'\in L\setminus\{0\}$ such that $y=x+x'$.
Let $f\in M\setminus\{0\}$ and $g,h\in C_b(U_f)_+$ be such that
$g$, $h$ and $h-g$ are representatives of $x$, $y$ and $x'$, respectively,
and both $g$ and $h$ have property~(C). We need to show that $h-g$ has property~(C)
as well.

Fix $t\in U_f\setminus\{t_0\}$. By assumption, there are $z=z_t\in N_f^*$,
$w=w_t\in N_f^*$ and $v\in N_f^*$ with the following properties
\begin{itemize}
\item[(i)]$z\leq g,\ w\leq h\,\text{ and }v\leq h-g\text{ on }U_f$;
\item[(ii)]$z=g,\ w=h\,$ on $V(t)$;
\item[(iii)]$z\ll g$, $w\ll h$ and $v\ll h-g$ on $V[t_0]$.
\end{itemize}
Note that $V(t)\cap V(t_0)=\emptyset$.

Since $g\le h$ on $U_f$, it follows from (ii) that $z\le w$
on~$V(t)$. Take $r\in G$, $r\geq0$ such that $r\geq w-z$ on
$V'(t)\subseteq V(t)$ and $r=0$ on ${}^cV(t)$. Then
\begin{equation*}
\begin{split}
&(w-z)\wedge r\in G,\quad  (w-z)\wedge r=w-z\,\text{ on }V'(t),\\
&(w-z)\wedge r\leq w-z\, \text{ on }U_f,\quad (w-z)\wedge
r\leq0\,\text{ on }{}^cV(t).
\end{split}
\end{equation*}
Put $d=0\vee((w-z)\wedge r)\in G$; then $0\leq d$, $d\leq w-z$ on~$V(t)$,
$d=w-z$ on $V'(t)$ and $d=0$ on ${}^cV(t)$. Let $e=d\vee v\in G$;
in fact, $e\in N_f$. Then $e\leq h-g$ on $U_f$, $e\ll h-g$ on $V[t_0]$
and $e=h-g$ on $V'(t)$. Thus $h-g$ has property~(C) and so $y-x\in C$.
\end{proof}

We next show that the functions with property (C) enjoy an even stronger
property.

\begin{lem}\label{approxlem}
Let\/ $h\in L_f$ have property\/ {\rm(C)}, where\/ $f\in M\setminus\{0\}$,
$f\le_M f_0$. For a compact subset\/ $K$ of\/ $X$ such that
$K\subseteq U_f$, there is\/ $z\in I_f(h)$ such
that\/ $z=h$ on~$K$. Moreover, if\/ $v\in I_f(h)$, then we can
choose\/ $z$ such that, in addition, $v\le_M z$.
\end{lem}

\begin{proof}
Since $K$ is compact, there exist $t_1,\ldots,t_n\in K$ such that,
for some $z_j=z_{t_j}\in N_f$, $z_j\leq h$ on $U_f$, $z_j=h$ on
$V(t_j)$, $z_j\ll h$ on $V[t_0]$ for all $1\leq j\leq n$ and
$K\subseteq V$, where $V=\bigcup_{j=1}^n V(t_j)$. Let $U_j$ be
open neighbourhoods of $t_0$ such that $0\ll z_j\ll h$ on
$U_j\setminus\{t_0\}$ and put $U=\bigcap_{j=1}^n U_j$. Note that
$U_j\cap V(t_j)=\emptyset\,$ so that $U\cap V=\emptyset$.

Letting $z'=z_1\vee \cdots \vee z_n \in N_f$ we have $z'\le h$ on
$U_f$, $z'=h$ on $K$ and $z'\ll h$ on $U\setminus\{t_0\}$. Take
$r\in M$ such that $r\geq h$ on some open neighbourhood $U'$ of
$t_0$ contained in $U$ and $r=0$ on~${}^cU$. Take $r'\in M$ such
that $r'\geq h$ on some open set $V'\subseteq V$ containing $K$
and $r'=0$ on ${}^cV$.
Letting $z=\frac12(z'\wedge r)+z'\wedge r'\in N_f$ we have $z\le h$
on $U_f$, $z=h$ on $K$ and
$z\ll z'\ll h$ on $U'\setminus\{t_0\}$. Thus $z\in I_f(h)$.

Now assume that $v\in I_f(h)$. Since both $z$ and $v$ belong to
the interval $I_f(h)$, there is $z''\in I_f(h)$ such that $z\le _M z''$
and $v\le _M z''$, and obviously $z''=h$ on $K$.
\end{proof}

The image of $\tau$ can now be identified.

\begin{prop}\label{image}
The monoid homomorphism\/ $\tau$ maps\/ $M'$ onto~$C$.
\end{prop}

\begin{proof}
We first show $\tau(M')\subseteq C$. Let $E\in\Lambda(N_f,D_f)\setminus\{0\}$,
where $f\in M\setminus\{0\}$, $f\leq_M f_0$. We have to verify that
$h=\tau_f(E)=\sup E$ has property~(C).

Let $t\in U_f$ and take $b\in N_f$ such that $b=1$ on $V(t)$ and
$b\leq_M1$; then $b\in D_f$. Let $E'\in\Lambda(N_f,D_f)$ and
$k\geq1$ be such that $E+E'=kD_f$. Take $z\in E$ and $z'\in E'$
with $z+z'=kb$. Letting $h'=\sup E'$ it follows that $h+h'=k$,
since $\tau_f$ is a monoid homomorphism. Consequently
\[
k=h+h'\geq z+z'=kb=k\quad\text{on }V(t)
\]
and so $z=h$ on~$V(t)$. Clearly $z\leq h$ on~$U_f$.

Let $r,\,r'\in M$ be such that $r\geq h$ on $V'(t)\subseteq V(t)$,
$r=0$ on ${}^cV(t)$ and $r'\geq h$ on $V'(t_0)\subseteq V(t_0)$,
$r'=0$ on ${}^cV(t_0)$ for some $V(t_0)$ such that
$V(t_0)\cap V(t)=\emptyset$. (This can be achieved by possibly shrinking the
neighbourhood~$V(t)$.) Letting $z_t=\frac12(z\wedge r')+z\wedge r$
we obtain $z_t\in N_f$ such that $z_t\leq h$ on $U_f$, $z_t=h$ on
$V'(t)$ and $z_t\ll h$ on $V'[t_0]$. Therefore, $h$ has
property~(C) and thus $\tau(\ol E)\in C$.

In order to establish the reverse inclusion, $C\subseteq\tau(M')$, note that
for each non-zero $y\in C$ there is $n\ge1$ such that $y<_L n1$, whence
$n1-y\in C$ by Lemma~\ref{oembedding}. There exists a
representative $h\in L_f$ of $y$, for some non-zero $f\in M$
with $f\le _M f_0$, such that both $h$ and $n1-h$ have property~(C).
Moreover, we can assume that $n1-h\gg\eps$ on $U_f$ for some $\eps>0$.
We will show that $I_f(h)+I_f(n1-h)=nD_f$ which, together with Property~(3)
above, entails that
$h=\sup I_f(h)=\tau_f(I_f(h))$ and thus $y\in \tau(M')$.

To this end take $0\ne g\in nD_f$.  Put $\delta =\eps/3$ and let
$K_1=\{t\in X\mid g(t)\le\delta/3\}$ and $K_2=\{t\in X\mid g(t)\ge\delta/2\}$.
For each rational number $\rho$ with
$0<n1-\rho<\delta$ we get, by using the assumption on the richness
in Urysohn functions, an element $a\in G$ such that $0\le a \le \rho$,
$a=\rho$ on $K_2$ and $a=0$ on $K_1$. Let $g'=g\wedge(n1-a)\in N_f$.
Note that $g'=g$ on $K_1$, $g'\le g$ and $g'\ll\delta$ on~$U_f$.
Set $K_3=\{t\in X\mid g(t)\ge\delta/3\}$. Then
$g-g'$ is supported on $K_3$. Take a non-zero function $v\in N_f$
supported on a compact neighbourhood $K_4$ of $t_0$ such that
$K_4\cap K_3=\emptyset $, and with $v \ll g' $ on
$K_4\setminus\{t_0\}$. Put $g''=g-g'+v$. Note that
$g'+v\ll2\delta\ll n1-h$ on~$U_f$. Set $h'=n1-h-g'$; then $v\in I_f(h')$
because $v+\delta\ll h'$ on $U_f$, and $h'$
has property~(C) by the arguments in Lemma~\ref{oembedding}.

Applying Lemma~\ref{approxlem} to $h$ and $h'$, respectively we get
$z_1\in I_f(h)$ and $z_2\in I_f(h')$ such
that $z_1=h$ on $K_3$, $z_2=h'$ on $K_3$ and $v \le _M z_2$. We
claim that $g''\le _M z_1+z_2$.
Indeed, on $K_3$, we have $z_1=h$ and $z_2=h'=n1-h-g'$,
so that $z_1+z_2=n1-g'\ge g-g'=g''$.
On ${}^cK_3$, we have $g''=v$, because $g-g'$
is supported on $K_3$ and $v\le_M z_2<_M z_1+z_2$. It follows that
$g''\le _M z_1+z_2$.

Since $M$ has the Riesz property, we get $g''=g_1+g_3$ with
$g_1,g_3\in N_f$ and $g_1\le _M z_1$ and $g_3\le _M z_2$. Set
$g_2=g_3+g'-v\in N_f$, and observe that
\[
g_2=g_3+g'-v\le _M z_2+g'-v\le z_2+g'\in I_f(n1-h)
\]
so that $g=g''+g'-v=g_1+g_2$ with $g_1\in I_f(h)$ and $g_2\in I_f(n1-h)$,
as desired.
\end{proof}

\begin{rema}\label{rem:image}
The proof of Proposition~\ref{image} in fact shows that, for every $h\in L_f$
with property~(C), the interval $I_f(h)$ is complemented.
\end{rema}

\section{Multiplier Algebras}\label{sec:multi-algebras}

In this section we will use the properties of the monoids studied above
to construct a \C* $A$ with the property that $\mloc\ne M_{\text{loc}}(\mloc)$.
Throughout let $(M,u)$ be a fixed monoid as considered in
Section~\ref{sec:monoids}. For every \C* $B$ of real rank zero, there is an
isomorphism between its lattice of closed ideals and the lattice of order-ideals of the
monoid $V(B)$ (\cite{El}, \cite[Theorem 2.3]{zhang}, \cite[Theorem 2.1]{perera}).

Let $A$ be the prime, unital \AF/ such that $(V(A),[1_A])=(M,u)$.
By the Blackadar--Handelman theorem, \cite[Theorem 6.9.1]{Black}, the trace simplex
of $A$ is precisely $M_1^+(X)$, the simplex of probability measures on~$X$.
By Lin's theorem \cite[Corollary 3.7]{lin}, $M(I)$ has real rank zero
for every closed ideal $I$ of~$A$.
For an element $f\in M$, $f\leq_M u$ denote by $I_f=\overline{ApA}$ the closed ideal
generated by a projection $p\in A$ such that $[p]=f$. The order-ideal $V(I_f)$
is precisely
$N_f=\{g\in M\mid\hbox{$g\leq_M nf$}\text{ for some } n\in\N\}$.
For a projection $q\in M(I_f)$ define
\[
\text{supp}(q)=\{t\in U_f\mid\tau_f (D)(t)\ne0\}\subseteq U_f,
\]
where $D$ is the interval in $\Lambda (N_f,D_f)$ corresponding to~$[q]$ 
via the canonical isomorphism between $V(M(I_f))$ and $\Lambda(N_f,D_f)$
\cite[Theorem~2.4]{perera} and
$\tau_f\colon \Lambda (N_f,D_f)\to L_f$ is the canonical map defined
in Section~\ref{sec:monoids}.

Consider a fundamental sequence $(I_{n})_{n\geq0}$ of closed ideals of $A$
(compare Section~\ref{sec:loc-multi-algebras}).
The closed ideals $I_n$ are assumed to be of the form $I_{f_n}$,
where $f_n\in M$ and $f_n\le_M f_{n-1}$ for $n\ge 1$,
and $f_0$ is the distinguished function in~$M$. The fact that
$(I_n)_{n\geq0}$ is a fundamental sequence means that, for each open
neighbourhood $V(t_0)$, there is some $n_0$ such that
$U_{f_n}\subseteq V(t_0)$ for all $n\ge n_0$. Given such a sequence
$(I_n)_{n\geq0}$ we have $\mloc =\dirlim M(I_n)$, with canonical
maps $\varphi_{m,n}\colon M(I_n)\to M(I_m)$, $n\leq m$ and
$\varphi_n\colon M(I_n)\to\mloc$, $n\ge0$.
Since the $\varphi_n$'s are isometric embeddings, we will subsequently
suppress them when no ambiguity can arise.

Let $I$ be the closed ideal of $B=\mloc$ generated by $p_0$,
where $p_0\in A$ is a projection such that $[p_0]=f_0$.
Then $V(I)=J$ is the unique minimal order-ideal of~$M'$,
cf.\ Proposition~\ref{prop:min-order-ideal}; thus $I$ is the unique minimal
closed ideal of $B$ (use Theorem~\ref{thm:from-old}) and so $M_{\text{loc}}(B)=M(I)$.
Our aim is to construct a sequence of
projections $(p_n)_{n\in\N}$ in $B$ such that $(p_n)_{n\in\N}$ is strictly
convergent in $M(I)$ to a projection $p\in M(I)$, but $p\notin B$.
This will ensure that $B\ne M_{\text{loc}}(B)$.

Put $B_n=M(I_n)\subseteq B$, $n\geq0$ and let $I_n'$ denote the closed
ideal of $B_n$ generated by~$I_0$.
Then $I_0^{}=I_0'\subseteq I_1'\subseteq\ldots\subseteq I_n'\subseteq\ldots\ $ and
$B_0\subseteq B_1\subseteq \ldots\subseteq B_n\subseteq\ldots\,$, and
\begin{equation}\label{eq:describe-B}
B=\ol{\bigcup\nolimits_{n=1}^{\,\infty}B_n}\quad\text{and}\quad
I=\ol{\bigcup\nolimits_{n=1}^{\,\infty}I_n'}.\tag{3.1}
\end{equation}
Moreover, note that $B_n=M(I_n')$ for all $n\geq0$ since $I_n$ is an ideal
in $I_n'$ and hence $B_n\subseteq M(I_n')\subseteq M(I_n)=B_n$.

\smallskip
The following easy observation enables us to manoeuver between different
multiplier algebras.

\begin{lem}\label{Iconvergence}
Let\/ $(x_n)$ be a sequence in the ideal\/ $I_0$ converging in the
strict\/ $I$-topology to\/ $x\in M(I)$. Then\/ $x\in M(I_0)$.
\end{lem}

\begin{proof}
Clearly $(x_n)$ converges in the strict $I_0$-topology too, with
limit $y\in M(I_0)$, say. Since $(x-y)I_0=0$, it suffices to show
that $cI_0\ne 0$ whenever $0\ne c\in M(I)$.

Let $c$ be a non-zero element in $M(I)$. There are a closed essential ideal $I'$
of $A$ such that $I'\subseteq I_0$ and $z\in M(I')\cap I$ such
that $\|cz\|=1$. Given $0<\eps<1/2$, there are another closed
essential ideal $I''\subseteq I'$ and $z'\in M(I'')$, $\|z'\|=1$ such that
$\|cz-z'\|<\eps $. Let $a\in I''$ with $\|a\|=1$ be such
that $\|z'a\|>1-\eps $. Then
\[
\|cza\|\ge \|z'a\|-\|cza-z'a\|>1-2\eps>0.
\]
Hence $c(za)\ne 0$, and $za\in M(I')I''\subseteq I''\subseteq I_0$ as desired.
\end{proof}

The next two results are at the core of our construction.

\begin{lem}\label{notinB}
Let\/ $(p_n)_{n\in\N}$ and\/ $(q_n)_{n\in\N}$ be increasing sequences
of projections in\/ $B$ such that, for each~$n$, $p_n,\,q_n\in B_n$ and\/
$p_n+q_n\in B_0$. Suppose that\/
$\hbox{\rm dist}(p_n,B_{n-1})\geq\delta$ for all\/ $n\ge1$ and some\/
$\delta>0$. If\/ $(p_n)_{n\in\N}$ and\/ $(q_n)_{n\in\N}$ converge
strictly in\/ $M(I)$ to projections\/ $p$ and\/ $q$, respectively
such that $p+q=1$ then\/ $p\notin B$ and\/ $q\notin B$.
\end{lem}

\begin{proof}
Suppose that $p\in B$. There is a projection $p'\in B_n$ for some
$n\in\N$ such that $\|p-p'\|<\delta$. Note that
\[
(p_{n+1}+q_{n+1})p=p_{n+1}
\]
and $(p_{n+1}+q_{n+1})p'\in B_0B_n=B_n$. Since
\[
\|p_{n+1}-(p_{n+1}+q_{n+1})p'\|=\|(p_{n+1}+q_{n+1})(p-p')\|
\leq\|p-p'\|<\delta,
\]
we obtain that $\text{dist}(p_{n+1},B_n)<\delta$ contradicting our hypothesis.
Therefore $p\notin B$ and so $q=1-p\notin B$.
\end{proof}

\begin{prop}\label{strictconverg}
With the same notation and caveats as above, take compact neighbourhoods\/
$K_n$, $n\geq0$ of\/ $t_0$ such that\/
$K_n\subseteq K_{n-1}\subseteq U_{f_{n-1}}\cup\{t_0\}$ for all\/ $n\geq1$.
Let\/ $(h_n)_{n\ge 0}$ be an increasing approximate identity consisting of
projections for~$I_0$. Let\/ $(h_n'')_{n\ge0}^{}$ be a sequence of
non-zero projections in\/ $I_0$ such that\/ $h_n''\le h_n-h_{n-1}$ and\/
$\text{\rm supp}(h_n'')\subseteq K_n$ for all~$n$ $($so that, in
particular, $h_n''\in I_n${}$)$. Set\/ $h_n'=h_0''+h_1''+\cdots+h_n''$.
Then\/ $(h_n')_{n\geq0}^{}$ converges in the strict topology of~$M(I)$.
\end{prop}

\begin{proof}
By identity~\eqref{eq:describe-B}, the fact that $I_0$ is an \AF/
and since $B_nI_0^{}B_n$ is dense in~$I_n'$,
every element in $I$ can be approximated by an element of the form
$\sum_{i=1}^ky_ie_iz_i$, for some projections $e_i\in I_0$ and
$y_i,z_i\in M(I_n)$. It suffices to consider the case $k=1$, hence assume that
$x=yez$ with $e\in I_0$ a projection and $y,z\in M(I_n)$. Given
$\eps>0$, we have to find $n_0$ such that, for $m>\ell\ge n_0$,
we have $\|(h_m'-h_{\ell}')x\|<\eps$.

For $m\ge n\ge1$, we have $h_m'-h_{n-1}'=\sum_{i=n}^m h_i''\in I_n$. The
sequence $(h'_m-h'_{n-1})_{m\ge n}$ converges in the strict
$I_0$-topology to an element $q\in M(I_0)$. This follows
because $h_i''\le h_i-h_{i-1}$ and $(h_i)_{i\geq0}$ is an approximate
identity for~$I_0$. Note that $\text{supp}(q)\subseteq K_n$
by hypothesis.

We note that $qy\in M(I_n)$ and that $M(I_n)$ is a \C* of real rank zero
(as $I_n$ is an \AF/).
Hence, given $0<\eta<1$, there is a projection $p\in M(I_n)qy$
such that $\|qy-qyp\|<\eta $. Then $p$ is equivalent to a
subprojection of $q$ and thus $\text{supp}(p)\subseteq K_n$.
It follows from Lemma~\ref{closure} below that $pe=pe'$ for some
projection $e' \in I_n$, so that $pe\in I_n$ and thus
\[
qypez\in M(I_n)I_nM(I_n)\subseteq I_n.
\]
From
$$\|qyez-qypez\|<\eta\, \|z\|$$
we find $\text{dist}(qyez, I_n)<\eta\,\|z\|$. Since this holds for
all $0<\eta<1$, it follows that $qyez\in I_n$.

For all $m\ge n$ we have $(h_m'-h_{n-1}')(yez)=(h_m'-h_{n-1}')q(yez)$
and, since $qyez\in I_n$ and $(h_k')_{k\geq0}^{}$ is a Cauchy sequence in the
strict topology of $M(I_0)$, there is some $n_0\geq n$ such that, for all
$m>\ell\ge n_0$,
$$\|(h_m'-h_\ell')yez\|=\|(h_m'-h_\ell')qyez\|<\eps,$$
as desired.
\end{proof}

The following somewhat technical lemma completes the proof of
Proposition~\ref{strictconverg}.

\begin{lem}\label{closure}
Let\/ $f\in M$ be such that\/ $f\le_M f_0$. Let\/
$p$ be a non-zero projection in\/ $M(I_f)$ such that the closure in\/
$X$ of\/ $\text{\rm supp}(p)$ is contained in\/ $U_f\cup\{t_0\}$. Then,
given a projection\/ $e$ in\/ $I_0$, there is a projection\/ $e'$ in\/
$I_f$ such that\/ $pe=pe'$. In particular, $pI_0\subseteq I_f$.
\end{lem}

\begin{proof}
Let $K$ denote the closure in $X$ of
$\text{supp}(p)$. By hypothesis, $K\subseteq U_f\cup\{t_0\}$.

Set $P=G\cap C(X)_+$; then $P$ is a countable dimension monoid.
Let $N$ be the order-ideal of
$P$ consisting of those functions $g\in P$ such that $g=0$ on~$K$.
Put $S=M+N\subseteq P$ and observe that $S$ is the
disjoint union of $N$ and $M\setminus \{0\}$. (Note that the
sum of an element in $N$ and a non-zero element in $M$ is
an element in $M$, because an element in $P$ is in $M$ if and
only if it is strictly positive on~$V[t_0]$.)

Since $S$ is a dimension monoid, there is an (up to
isomorphism unique) unital \AF/ $D$ such that $(V(D),[1_D])\cong (S,u)$.
We have a monoid homomorphism\break
$\lambda\colon(M,u)\to(S,u)$ and therefore there is a unital *-homomorphism
$\psi\colon A\to D$ with the property $V(\psi)=\lambda$.
Let $I'$ be the closed ideal of $D$ such that $V(I')=N$.

Since $p\in M(I_f)$, $p$ is the limit with
respect to the strict topology of $I_f$ of an increasing sequence $(h_n)$
of projections, which forms an approximate identity for
$pI_fp$. Let $I$ denote the closed ideal of $A$ generated by $(h_n)$, that is,
$I=\overline{\bigcup_{n=1}^{\infty}Ah_nA}$.
Then $I\subseteq I_f$ and clearly $p$ is also the strict limit of
$(h_n)$ with respect to~$I$. Observe that
$\text{supp}(h_n)\subseteq K$ for all~$n$, so that $\lambda(V(I))$ is
an order-ideal of~$S$. Let $I''$ be the closed ideal of
$D$ corresponding to $\lambda(V(I))$, that is, the closed ideal
generated by all the projections $q$ in $D$ such that
$[q]\in\lambda (V(I))$. The map $\rest\psi I\colon I\to I''$
induces an isomorphism of monoids $V(\rest\psi I)\colon V(I)\to V(I'')$
sending the canonical interval in $V(I)$ onto the
canonical interval in $V(I'')$. Since both $I$ and $I''$ are
\AF/s, it follows from Elliott's theorem that $\rest\psi I$ is an
isomorphism from $I$ onto~$I''$. Since $A$ is prime and $I$ is a
non-zero ideal of $A$, we conclude that the map $\psi\colon A\to D$
is injective. Therefore we can identify $A$ with a
\Cs* of $D$ via~$\psi$. Note that, under this
identification, $I$ is a closed ideal of $D$ such that $II'=0$.
Thus $I'\subseteq I^{\perp}$, where we denote by $I^{\perp}$ the
orthogonal ideal of $I$ in~$D$, that is, the set of all elements
in $D$ that annihilate~$I$.

Now we need a suitable decomposition for the projection $e\in I_0$.
Let $f'$ be a fixed non-zero element of $V(I)$, so that
$U_{f'}\subseteq K$. There is a compact neighbourhood $K_1$ of
$t_0$ such that $K_1\subseteq U_{f'}\cup \{t_0\}$. Let $W$ be an
open neighbourhood of $t_0$ such that $W\subseteq K_1$ and
$[e]=\beta f_0$ on~$W$, for some $\beta\in\Q_+\setminus\{0\}$.
We can select a compact subset $K_2$ of $X$ and an open
subset $U$ of $X$ such that
$K\subseteq U\subseteq K_2\subseteq U_f\cup \{t_0\}$,
so that we have the following situation:
\[
t_0\in W\subseteq K_1\subseteq U_{f'}\cup \{t_0\}\subseteq
K\subseteq U\subseteq K_2\subseteq U_f\cup \{t_0\}.
\]
Using suitable Urysohn functions, we will establish an orthogonal decomposition
$e=e_1+e_2$, where $e_1,e_2$ are projections in $A$ with $e_1\in
I_f$ and $\text{supp}([e_2])\subseteq K_1\cup K'$, where $K'$ is a
compact subset of $X$ such that $K\cap K'=\emptyset$.

Put $g=[e]\in M$. Take $\rho\in\Q$ such that $\rho\gg g$ on~$X$.
Then there is $r_1\in M$ and an open subset $V$
of $X$ with $K\subseteq V\subseteq U$ such that $r_1= \rho$ on
$V$ and $r_1=0$ on $X\setminus U$. Set $g_1'=g\wedge r_1$
and note that $g_1'\in M$ and $\text{supp}(g_1')\subseteq
K_2\subseteq U_f\cup \{t_0\}$. It follows
that there is a positive integer $k$ such that $g_1'\le _M kf$ and
hence $g_1'\in N_f$.

Take a rational number $\alpha$ such that $0<\alpha <\beta$,
so that $\alpha f_0\ll g$ on $W\setminus\{t_0\}$. In addition take
$\rho'\in\Q$ such that $\alpha f_0\ll\rho'$ on~$W$.
There is a Urysohn function $r_2\in M$ such that $r_2=\rho'$ on an
open neighbourhood $W'\subseteq W$ of $t_0$, and $r_2=0$ on $X\setminus W$.
Set $g_2'=\alpha f_0\wedge r_2\in M$
and note that $g_2'\le_M g_1'$ and $\text{supp}(g_2')\subseteq K_1$.
It follows that there is $\ell\in\N$ such
that $g_2'\le_M \ell f'$, so that $g_2'\in V(I)$. Now consider
the element $g_2=g_2'+(g-g_1')$. Note that $g_2\in P$ and $g_2\gg0$
on $W'\setminus \{t_0\}$ so that $g_2\in M$. Finally we set
$g_1=g_1'-g_2'\in M$. Then $g=g_1+g_2$ with $g_1,g_2\in M$, and we
have $g_1\in N_f$ and $\text{supp}(g_2)\subseteq K_1\cup K'$,
where $K'=X\setminus V\subseteq X\setminus K$. There is a
corresponding orthogonal decomposition $e=e_1+e_2$, with $e_1\in I_f$
and $[e_2]=g_2$.

The element $g_2\in M$ decomposes in $S=M+N$ as
$g_2=g_2'+(g-g_1')$, where $g_2'\in N_{f'}$ and $g-g_1'\in N$,
because $g-g_1'$ vanishes on~$K$. This implies an orthogonal
decomposition $e_2=e_2'+e_2''$ of $e_2$ in $D$ such that
$[e_2']=g_2'$ and $[e_2'']=g-g_1'$ in $V(D)=S$. Since $f'\in V(I)$,
we have $N_{f'}\subseteq V(I)$, and we know that the closed ideal
of $D$ generated by its order-ideal $V(I)$ is precisely $I$, so
that $e_2'\in I $. On the other hand, $e_2''\in I'$ because
$[e_2'']\in N=V(I')$. Therefore we can write
\[
e=e_1+e_2=(e_1^{}+e_2')+e_2'',
\]
where $e_1+e_2'\in I_f$ and $e_2''\in I'$. Set $e'=e_1^{}+e_2'\in I_f$
and $e''=e_2''\in I'$.

Note that $I\oplus I^ {\perp}$ is an essential closed ideal of $D$ and we
have a canonical inclusion $\iota \colon D\to M(I\oplus
I^{\perp})=M(I)\oplus M(I^{\perp})$. The sequence
$(h_n,0)$ converges in the strict topology of $I\oplus I^{\perp}$
to $(p,0)\in M(I)\oplus M(I^{\perp})$. We have the following commutative diagram
\[
\xymatrix{D\ar[r]^(.28){\iota} &\ \; M(I)\oplus M(I^{\perp}) \ar[d]^{\pi_1}\\
\ar[u]^{\psi} A \ar[r] & {M(I)}}
\]
where $\pi_1\colon M(I)\oplus
M(I^{\perp})\to M(I) $ is the projection onto the first component.

Using the above decomposition $e=e'+e''$ and
$I'\subseteq I^{\perp}$, we find that the image of $e$ in $M(I)\oplus M(I^{\perp})$
is $(\rest{e'}{I},\rest{e}{I^{\perp}})$, where, for a closed ideal
$J$ of $D$, we denote by $\rest xJ$ the image of $x$ in $M(J)$ under
the canonical restriction map $D\to M(J)$. It follows that
\[
pe=\pi_1\bigl((p,0)(\rest{e'}{I},\rest{e}{I^{\perp}})\bigr)
=\pi_1\bigl((\rest{pe'}{I},0)\bigr)=pe'\in M(I).
\]
Since $p\in M(I_f)$ and $e'\in I_f$ we get $pe=pe'\in I_f$, as claimed.
\end{proof}

Although not strictly necessary, we shall assume for the remainder of this paper
that $X=[0,1]$, $t_0=0$ and $f_0(t)=t$, $t\in X$. This allows us to construct certain
discontinuous functions without undue notational complications.

In the following lemma we denote by
$\varphi_{f_1,f_2}\colon M(I_{f_2})\to M(I_{f_1})$, $f_1\le _M f_2$
the canonical restriction map.

\begin{lem}\label{monsterlemma}
Let\/ $A$ be the \AF/ such that\/ $(V(A),[1_A])=(M,u)$.
Let\/ $f$ be a non-zero element in\/ $M$ with\/ $f\le_M f_0$, and
let\/ $p$ be a non-zero projection in~$I_0$. Then there
exist\/ $f'\le_Mf$ and a projection\/ $q\in M(I_{f'})$ with\/
$q\le\varphi_{f',f_0}(p)$ such that\/ $q$ is not equivalent in\/
$M(I_{f'})$ to any projection in\/ $\varphi_{f',f}(M(I_f))$.
\end{lem}

\begin{proof}
We can assume that $f=\lambda f_0$ and $h=2\lambda' f_0$ for some
positive rational numbers $\lambda$ and $\lambda'$, where $h=[p]$.
Let $f'\in M$ be a non-zero function such that $f'\le _M f$ and
$U_{f'}\varsubsetneq U_f=(0,1]$. Let $\rho$ be the left-most positive
number $t$ such that $f'(t)=0$. Let $0<\mu\in\Q$ be such that
$\mu<\lambda'$ and choose a strictly  increasing sequence of positive rational
numbers $(\mu_n)_{n\ge 1}$ converging to~$\rho$. Define a function
$g$ on $U_{f'}$ as follows. On $U_{f'}\setminus(0,\rho)$ we set
$g=0$. In the interval $(0,\mu_1]$, we
set $g(t)=\mu t$, so that $g(\mu_1)=\mu\mu_1\in\Q$. In
the interval $[\mu _1,\mu _2]$, define $g$ to be the restriction
to $[\mu _1,\mu _2]$ of a Urysohn function $r_1\in G$ such
that $r_1=\mu_1\mu$ on $[0,\mu_1]$, $r_1=0$ on $[\mu_2,1]$
and $0\le r_1\le \mu_1\mu$. Observe that $g=(\mu f_0)\wedge r_1\in M$
on $[0,\mu_2]$. On $[\mu_2,\mu_3]$,
we define $g$ as the restriction to $[\mu_2,\mu_3]$ of a Urysohn
function $r_2\in G$ such that $r_2=0$ on $[0,\mu _2]$,
$r_2= \mu _1\mu $ on $[\mu _3,1]$ and $0\le r_2\le \mu
_1\mu$. Note that $g=((\mu f_0)\wedge r_1)\vee r_2\in M$
on $[0,\mu _3]$. We continue in this way, obtaining a continuous
function $g$ on $(0,\rho)$ which cannot be extended to a continuous
function at $\rho$. Observe that, by the above arguments, the
function $g$ is locally in $M$, that is, for each $t\in U_{f'}$
there is a function $z_t\in M$ such that $g=z_t$ on an open
neighbourhood of~$t$. It follows easily that $g$ has
property~(C), and by construction $g\ll h_1$ on $U_{f'}$,
where $h_1=\frac12 h$.

Let $E_1=I_{f'}(g)$ and $E_2=I_{f'}(h_1-g)$. By Lemmas~\ref{lem:ifh-interval}
and~\ref{oembedding} together with Remark~\ref{rem:image},
$E_1,E_2\in\Lambda(N_{f'},D_{f'})$ and $\tau_{f'}(E_1)=g$ and
$\tau_{f'}(E_2)=h_1-g$. It follows that
$E=E_1+E_2\in\Lambda(N_{f'},D_{f'})$ and
$\tau_{f'}(E)=\tau_{f'}(E_1)+\tau_{f'}(E_2)=h_1$. We claim that
$E+([0,h_1]\cap N_{f'})=[0,h]\cap N_{f'}$. Let $z\in N_{f'}$ be such that
$z\le_M h=2h_1$. By the Riesz property
of $M$, we can write $z=z_1+z_2$ with $z_1\le _M h_1$ and $z_2\le_M h_1$.
Since $h_1=\lambda 'f_0$ with $\lambda ' >0$, and
$f'(\rho)=0$, we conclude that $z_1<_M h_1$ and $z_2<_M h_1$. Take
$z_1'\in M$ such that its support is contained in a
closed interval of the form $[0,\beta]$, with $\beta <\rho$, and
$z_1'<_M z_1$ and $z_1-z_1'<_M h_1-z_2$. (This is possible because
there is $\eps>0$ such that $z_1\ll\eps$ and $h_1-z_2\gg\eps$
on $[\rho-\beta, 1]$ for some $\beta <\rho$.) Since $E$
contains all functions in $[0,h_1]\cap N_{f'}$ whose support is
contained in a closed interval $[0,\beta]$ with $\beta<\rho$
(because of the special construction of~$g$), it follows that
$z_1'\in E$. On the other hand,
$z_2+(z_1-z_1')<_M z_2+(h_1-z_2)=h_1$ and so
\[
z=z_1+z_2=z_1'+(z_2+(z_1-z_1'))\in E+([0,h_1]\cap N_{f'}).
\]
This shows that $E+([0,h_1]\cap N_{f'})\supseteq[0,h]\cap N_{f'}$ and the
reverse inclusion is obvious.

Since
\[
E_1+(E_2+([0,h_1]\cap N_{f'}))=[0,h]\cap N_{f'},
\]
we find that there is a projection $q\in M(I_{f'})$ such that
$q\le\varphi_{f',f_0}(p)$ and the interval in $\Lambda(N_{f'},D_{f'})$
corresponding to $[q]$ is~$E_1$. Since
\[
\tau _{f'}(E_1)=\tau _{f'}(I_{f'}(g))=g,
\]
and $g$ cannot be extended to a continuous function on $U_f$, we
infer from the commutative diagram~(2.1)
that $q$ is not equivalent in $M(I_{f'})$ to a projection
in $\varphi_{f',f}(M(I_f))$.
\end{proof}

With the same notation as in Proposition~\ref{strictconverg}, let
$p$ be the strict limit of the sequence $(h_n')$. By Lemma~\ref{Iconvergence},
$p\in M(I_0)=B_0$. We will now put all the above ingredients together
to obtain sequences of projections satisfying the requirements in
Lemma~\ref{notinB}.

\begin{prop}\label{monsters}
With the same notation and caveats as above, set\/ $p'_n=p-h_{n-1}'\in B_0$
for\/ $n\geq1$ and $p'_0=p$. Then there are\/ $\delta>0$, a sequence\/ $(f_n')$
of elements of\/ $M$ with\/ $f_0'=f_0^{}$ and\/ 
$f_n'\le_M f_{n-1}'\le _M f_{n-1}$ for all\/ $n\geq1$ as well as
orthogonal decompositions\/ $1=p'_n+p_n+q_n$ such that\/ $p_n,q_n\in M(I_{f'_n})$
for all\/ $n\ge0$ and\/ 
{\rm$\text{dist}(p_n,M(I_{f'_{n-1}}))\geq\delta$} for all\/ $n\ge1$.

It follows that the sequences\/ $(p_n)$ and\/ $(q_n)$ converge
in the strict $I$-topology to\/ $e\in M(I)$ and\/ $1-e$, respectively.
\end{prop}

\begin{proof}
Let $\delta>0$ be such that, for all \C*s $C,\,D$ with $C\subseteq D$ and for all
projections $e'\in D$, $\text{dist}(e',C)<\delta$ implies that $e'$
is equivalent to a projection in~$C$; cf. \cite[Lemma~6.3.1
and Proposition~2.2.4]{rordam}.

The sequences $(p_n)$ and $(q_n)$ are constructed inductively. To
start with we set $p_0=1-p$, $q_0=0$. Then $p_0,q_0\in B_0$ and
$1=p'_0+p_0+q_0$. Suppose that, for $n\ge0$, we have an orthogonal
decomposition $1=p'_n+p_n+q_n$ satisfying the stated conditions.
Note that $p'_n=p'_{n+1}+h_n''$.

By Lemma \ref{monsterlemma}, there exist a non-zero $f'_{n+1}\in
M$ with $f'_{n+1}\le _M f'_n$ and $f'_{n+1}\le _M f_{n+1}$, and a
projection $q'_{n+1}\in M(I_{f_{n+1}'})$ such that
$q'_{n+1}\le\varphi_{f'_{n+1},f_0^{}}(h''_n)$ in $M(I_{f'_{n+1}})$ and $q'_{n+1}$
is not equivalent in $M(I_{f'_{n+1}})$ to a projection in
$\varphi_{f'_{n+1},f'_n}(M(I_{f'_n}))$. The latter condition implies that
$\text{dist}(q'_{n+1},M(I_{f'_n}))\geq\delta$.

Identify all the projections constructed so far with their images
in $M(I_{f'_{n+1}})$ under the canonical inclusions; then
$q'_{n+1}\le h_n''$. Set $p_{n+1}=p_n+q'_{n+1}$ and
$q_{n+1}=q_n+(h''_n-q'_{n+1})$. Then
\[
\text{dist}(p_{n+1},M(I_{f'_n}))=\text{dist}(q'_{n+1},M(I_{f'_n}))\geq\delta;
\]
moreover,
\begin{align*}
1&=p'_n+p_n^{}+q_n^{}=p'_{n+1}+h''_n+p_n^{}+q_n^{}\\
 &=p'_{n+1}+q'_{n+1}+(h_n''-q'_{n+1})+p_n^{}+q_n^{}\\
 &=p'_{n+1}+p_{n+1}^{}+q_{n+1}^{}.
\end{align*}
This concludes the inductive construction.

Now we consider all the projections as projections in $M(I)$ and all the
algebras $M(I_{f'_n})$ as \Cs*s of~$M(I)$. It is a simple
matter to show that $(p_n)$ converges in the strict
$I$-topology. Indeed, fix $a\in I$ and $\eps>0$. Then there is $n_0$
such that $\|p'_{n_0}a\|=\|(p-h'_{n_0-1})a\|<\eps$. For $m>m'\ge n_0$, we have
$p_m-p_{m'}\le p'_{n_0}$ and so $\|(p_m-p_{m'})a\|<\eps$.
Similarly, the sequence $(q_n)$ is strictly convergent. Let $e$
be the strict limit of $(p_n)$. Since $p_n+q_n=1-p'_n$ converges strictly to $1$,
it follows that $(q_n)$ converges strictly to $1-e$.
\end{proof}

We are ready to complete the proof of our main result, Theorem~\ref{thm:main}.

\medskip\noindent
By Proposition~\ref{monsters}, we can construct projections $p,\,q$
in $M(I)$ such that $p+q=1$ and $p$ and $q$ are strict limits of
sequences $(p_n)$ and $(q_n)$, respectively satisfying the conditions stated
in Lemma~\ref{notinB} with respect to the \Cs*s
$B_n=M(I_{{f'_n}})$. Since $(I_{f'_n})$ is a fundamental
sequence of closed ideals of $A$ with $I_{f'_0}=I_0$, it follows from Lemma~\ref{notinB}
that $p,\,q\notin B$. Therefore
\[
M_{\text{loc}}(\mloc )=M_{\text{loc}}(B)=M(I)\ne B=\mloc.
\]

\section{Local Multiplier Algebras}\label{sec:loc-multi-algebras}

In this section we add a few remarks on a systematic approach to 
understanding the ideal structure of $\mloc$. Let $A$ be a separable prime \C*. 
Then $A$ is primitive \cite[4.3.6]{pedbook} and hence $0\in\Prim$. Since
$\Prim$ is second countable \cite[4.3.4]{pedbook}, we can find a countable basis
$(U_n)$ of open neighbourhoods of~$0$, with $U_{n+1}\subseteq U_n$ for all~$n$. 
These open sets correspond to a cofinal
countable family $(I_n)$ of non-zero closed ideals of $A$ such that
$I_{n+1}\subseteq I_n$ for all~$n$. We call such a sequence
$(I_n)$ a {\it fundamental sequence of ideals\/} of~$A$. Obviously,
we have $\mloc=\dirlim M(I_n)$ for such a fundamental sequence
of ideals.

The sequence $(I_n)$ determines a fundamental sequence $(J_n)$ in
$\mloc$, so that $0$ has a countable basis of
neighbourhoods in $\text{Prim}(\mloc)$. Indeed, define $J_n$ as
the closed ideal of $\mloc $ generated by~$I_n$. Given a non-zero
closed ideal $J$ of $\mloc $ we obtain that $J\cap A$ is a
non-zero ideal of~$A$ \cite[2.3.2]{AM2}; hence there is $m$ such that
$I_m\subseteq J$, and so $J_m\subseteq J$. It follows that
$\Mloc{\mloc}=\dirlim M(J_n)$. In general, the iterated local multiplier 
algebra $M_{\rm loc}^{(k)}(A)$ can be computed by taking the direct limit
$\dirlim M(I_n^{(k-1)})$, where $I_n^{(k-1)}$ is the closed
ideal of $M_{\rm loc}^{(k-1)}(A)$ generated by~$I_n$.

We can distinguish three different types of behaviour. The first one
corresponds to the case where all ideals $J_n$ are equal to
$\mloc$, that is, $\mloc$ is a simple \C*. Of course, this
happens when $A$ is simple and unital, but it can also occur when
all the ideals $I_n$ are different; examples of this behaviour were 
constructed in~\cite{AM1}. 
A second possibility is that $0$ is an isolated point
in $\text{Prim}(\mloc)$ which has more than one point. This
is the case if and only if the sequence $(J_n)$ stabilizes and 
$J_n\ne\mloc$ for large~$n$. In this case we have 
\[
\Mloc{\mloc}=M(J_{n_0}),
\]
where $J_{n_0}=J_n$ for all $n\ge n_0$, and so
$M_{\rm loc}^{(k)}(A)=\Mloc\mloc$ for all $k\ge2$.
Our examples in the present paper are of this type.

The third kind of prime \C*s consists of those such that the
family $(J_n)$ is strictly decreasing. Although we will not go
into the details of the construction, it is possible to give
explicit examples of \AF/s in this class using the methods
developed in the present paper. However it seems technically
challenging to analyze the possible lack of stabilization of the
increasing chain $\bigl(M_{\rm loc}^{(k)}(A)\bigr)_{k\in\N}$.

\end{document}